\documentstyle[12pt]{amsart}
\newtheorem{theorem}{Theorem}

\newtheorem{lemma}{Lemma}
\newtheorem{corollary}{Corollary}
\title[Natural Transformations]
{Natural Transformations\\
from Constructible Functions to Homology}
\author[G. Kennedy]{Gary Kennedy}
\address{Ohio State University at Mansfield\\
Mansfield, Ohio 44906}
\email{kennedy@@math.ohio-state.edu}
\author[C. McCrory]{Clint McCrory}
\address{Department of Mathematics\\
University of Georgia\\Athens, Georgia 30602}
\email{clint@@math.uga.edu}
\author[S. Yokura]{Shoji Yokura}
\address{Department of Mathematics\\
College of Liberal Arts and Sciences\\
University of Kagoshima\\
1-21-30 Korimoto, Kagoshima 890\\
Japan}
\email{f77446@@sinet.ad.jp}
\thanks{The second author was partially supported by
NSF grant DMS-9403887.
The third author was partially supported by
Grant-in-Aid for Scientific Research No.\ 06640162,
Japanese Ministry of Education, Science and Culture.}
\thanks{We thank Robert Varley for many helpful comments.}
\date{\today}
\newcommand{\mpc}{\bold{c}} 
\newcommand{\cf}{\bold{1}} 
\newcommand{\tq}{\otimes\bold{Q}} 
\begin{document}
\maketitle
\begin{abstract}
For complex projective varieties, all natural transformations from
constructible functions to homology (modulo torsion) are linear
combinations of the MacPherson-Schwartz-Chern classes. In particular, the total
Chern class is the only such natural transformation $\mpc$ such that
for all projective spaces $\bold{P}$, the top component of
$\mpc(\cf_{\bold{P}})$ is the
fundamental class of $\bold{P}$.
\end{abstract}
%
%
\section{Introduction}
\label{intro}
\par
The MacPherson-Schwartz-Chern class natural
transformation $\mpc$ from the constructible functions
functor to homology \cite{macp} \cite{bs} satisfies
a remarkably stringent
normalization requirement:
for each nonsingular variety $X$, the element
$\mpc(\cf_X)$ of homology
assigned to the characteristic
function of $X$ is the total Chern class
$c(TX) \frown [X]$.
Now suppose that we entirely abandon this requirement.
Then each individual component $\mpc_i$
of the Chern class natural transformation
(assigning to the characteristic function
$\cf_X$ of a nonsingular variety the
homology class $c_i(TX) \frown [X]$)
is likewise a natural transformation,
as is any linear combination of these components.
\par
We will show that, modulo torsion, these linear combinations
are the only natural transformations between the two functors.
In particular, the MacPherson-Schwartz-Chern class is the only such
natural transformation satisfying
this weak normalization requirement:
for each projective space $\bold{P}$,
the top-dimensional component of
$\mpc(\cf_{\bold{P}})$ is the fundamental class $[\bold{P}]$.
We conjecture that the same statements are valid even
for integral homology.
We want to note two features of our proofs.
First, we never appeal to resolution of singularities.
Second, several of our arguments are similar to those
in the proof (attributed to A. Landman)
of the last Proposition of \cite{fulton}.
\par
%
%
\section{The functors}
\label{func}
\par
Consider the category of complex projective algebraic varieties.
If $X$ is such a variety, the
{\em characteristic function} of a subvariety $Z$
is the function $\cf_Z$ on $X$
whose value on $Z$ is 1
and whose value elsewhere is 0.
If $f: X \to Y$ is a morphism, then
$f_* \cf_Z$
is the constructible function
on $Y$ whose value at $y$
is the Euler characteristic of
$f^{-1}(y) \cap Z$
(a subvariety of $X$).
A finite linear combination (over $\bold{Z}$) of characteristic
functions of subvarieties of $X$ is called a {\em constructible}
function on $X$.
We define the pushforward of an arbitrary constructible
function by extending the previous definition by linearity,
and thus define the {\em constructible functions functor}
$\cal C$ to the category of abelian groups.
\par
We also want to work with an appropriate homology theory.
This will be either of the following:
\begin{itemize}
\item {\bf Ordinary singular or simplicial homology.}
The fundamental class of an $n$-dimensional variety $X$
is an element of  $H_{2n}(X)$.
\item {\bf Algebraic cycles modulo rational equivalence.}
The standard reference is \cite[Ch.\ 1]{int}.
The fundamental class of an $n$-dimensional variety $X$
is an element of the cycle class group
$A_n(X)$.
\end{itemize}
\par
%
%
\section{Projective spaces}
\label{proj}
\par
Suppose that  $\tau$  is a natural transformation
from $\cal C$
to  $H_{2i}$ or to $A_i$.
Suppose that  $n > i$;
then  $\bold{P}^i$  is naturally a linear subspace of  $\bold{P}^n$.
\begin{theorem}
\label{proj-re}
$\tau(\cf_{\bold{P}^n}) =
\binom{n+1}{i+1}\tau(\cf_{\bold{P}^i})$.
\end{theorem}
\begin{pf}
Consider the morphism $f:\bold{P}^n \to \bold{P}^n$,
\begin{equation*}
f: [x_0,x_1,\dots,x_n] \mapsto [x_0^2,x_1^2,\dots,x_n^2].
\end{equation*}
Let
\begin{equation*}
L_k = \{[x_0,x_1,\dots,x_n] \mid
\text{ at least } n-k \text{ of the coordinates equal zero}\}
\end{equation*}
and
\begin{equation*}
U_k = \{[x_0,x_1,\dots,x_n] \mid
\text{ exactly } n-k \text{ of the coordinates equal zero}\}.
\end{equation*}
Then $L_k$ is a union of $\binom{n+1}{k+1}$
linear subspaces of dimension $k$, and $U_k = L_k \setminus L_{k-1}$.
The restriction of $f$ to $U_k$
is a self-covering map of degree $2^k$.
Hence the restriction of $f$ to each component of $L_k$
is a branched self-cover of the same degree,
and the induced map
$f_*:H_{2k}(\bold{P}^n) \to H_{2k}(\bold{P}^n)$,
or the induced map
$f_*:A_k(\bold{P}^n) \to A_k(\bold{P}^n)$,
is multiplication by $2^k$.
Consider the constructible function $\cf_{U_k}=\cf_{L_k}-\cf_{L_{k-1}}$.
By naturality of $\tau$ we have
\begin{equation*}
2^k \tau(\cf_{U_k})
= \tau(2^k \cf_{U_k})
= \tau f_*(\cf_{U_k})
=  f_* \tau(\cf_{U_k})
= 2^i \tau(\cf_{U_k}).
\end{equation*}
If  $i \neq k$,
then $\tau(\cf_{U_k}) = 0$.
Therefore
\begin{equation*}
\tau(\cf_{\bold{P}^n})
=\tau\left(\cf_{L_i}+\sum_{k=i+1}^n \cf_{U_k}\right)
= \tau(\cf_{L_i}).
\end{equation*}
And
\begin{equation*}
\tau(\cf_{L_i})=\tau(\sum_K \cf_K),
\end{equation*}
where the sum is taken over the components of
$L_i$, since the two functions in question
differ only on a variety of dimension $i-1$.
Since $L_i$ is a union of $\binom{n+1}{i+1}$
linear subspaces of dimension $i$,
the statement of the theorem follows.
\end{pf}
\par
%
%
\section{Galois coverings}
\label{gct}
\par
A morphism $X \to Y$ of (irreducible) varieties is called {\em Galois}
if it is finite and surjective, and the
corresponding field extension $K(Y) \subset K(X)$ is Galois.
\begin{lemma}
\label{gct1}
Suppose that $X$ is a projective variety of dimension $n$.
Then there exists a normal projective variety $Z$,
a finite morphism $Z \to X$,
and a finite morphism $X \to \bold{P}^n$,
such that the composition $Z \to \bold{P}^n$ is Galois.
\end{lemma}
\begin{pf}
To construct a finite morphism $X \to \bold{P}^n$,
first embed $X$ in a projective space $\bold{P}^N$
and then project it from a general $\bold{P}^{N-n-1}$.
The field extension $K(\bold{P}^n) \subset K(X)$
is finite; hence it can be obtained by adjoining
a single element satisfying a minimal polynomial $p$.
Let $\cal K$ be the splitting field of $p$.
Then $\cal K$ is a normal extension of
$K(\bold{P}^n)$ \cite[Thm.\ 8.4, p.\ 82]{stewart}.
Let $Z$ be the normalization of $X$ in $\cal K$
\cite[III.8, Thm.\ 3]{mumford}.
Then $Z$ is projective
\cite[III.8, Thm.\ 4]{mumford}.
Since $K(Z) = \cal K$, the morphism $Z \to P^n$ is Galois.
\end{pf}
\par
\begin{lemma}
\label{gct2}
If $Z$ and $Y$ are normal varieties and $\gamma:Z \to Y$ is Galois,
with Galois group $G$, then $Y$ is isomorphic to $Z/G$.
\end{lemma}
\begin{pf}
The finite group $G$ acts by birational maps.
Let $g:Z \to Z$ be any one of these maps;
let $W \subset Z \times Z$ be the closure of its graph.
The projection $W \to Z$ onto the first factor is
a birational morphism with finite fibers.
By Zariski's Main Theorem
\cite[Ch.\ 3, sec.\ 9]{mumford}
this projection is an isomorphism.
Hence $g$ is in fact a morphism.
\par
The morphisms $\gamma$ and $\gamma \circ g$
agree on an open dense subset of $Z$; hence they
are equal.
The quotient $Z/G$ is a projective variety
\cite[pp.\ 126--127]{harris}.
By the universal property of $Z/G$,
there is a morphism to $Y$.
Again by Zariski, this morphism is an isomorphism.
\end{pf}
\par
\begin{lemma}
\label{gct3}
Under the same hypotheses,
$\gamma_*$ maps the invariant
subgroup $(H_*Z \tq)^G$
isomorphically to $H_*(Z/G)\tq$.
Likewise $\gamma_*$ maps
$(A_*Z \tq)^G$ isomorphically to $A_*(Z/G)\tq$.
\end{lemma}
\begin{pf}
Triangulate the quotient map \cite{hardt}.
The $G$-invariant simplicial homology of $Z$
is isomorphic to the homology
of the complex of $G$-invariant simplicial chains of $Z$.
(The proof is by averaging over the group $G$.)
And the $G$-invariant simplicial
chain complex of $Z$ is isomorphic to
the simplicial chain complex of $Z/G$.
\par
For the statement about rational equivalence groups,
see \cite[Example 1.7.6]{int}.
\end{pf}
\par
%
%
\section{The natural transformations}
\label{nt}
\par
\begin{theorem}
\label{main}
Suppose that $\sigma$ is a natural transformation
from the constructible function functor $\cal C$
to $H_* \tq$, ordinary singular homology with rational
coefficients, or to $A_* \tq$, rational equivalence theory
with rational coefficients.
Suppose that for each projective space $\bold{P}$
the top-dimensional component of
$\sigma_{\bold{P}}(\cf_{\bold{P}})$ vanishes.
Then $\sigma$ is identically zero.
\end{theorem}
\begin{corollary}
Suppose that $\tau$ is a natural transformation
from $\cal C$
to $H_* \tq$ or to $A_* \tq$.
Then $\tau$ is a rational linear combination
$\sum_{i=0}^{\infty} r_i \mpc_i$
of the components of the MacPherson-Schwartz-Chern class.
\end{corollary}
\par
Note that although the sum appearing in the corollary is
nominally infinite, for any particular variety it is a finite sum.
\par
\begin{pf*}{Proof of corollary}
Define $r_i$ to be the rational number for which
\begin{equation*}
\tau\left(\cf_{\bold{P}^i}\right)
= r_i \lbrack\bold{P}^i\rbrack + \text{ terms of lower dimension}.
\end{equation*}
Then apply the theorem to
\begin{equation*}
\sigma = \tau - \sum_{i=0}^{\infty} r_i \mpc_i.
\end{equation*}
\renewcommand{\qed}{}
\end{pf*}
\begin{pf*}{Proof of theorem}
Suppose that $Z$ is a subvariety of $X$.
If $\sigma(\cf_Z)$ is zero in the homology group
of $Z$, then by naturality it is likewise zero in the
homology group of $X$.
Hence it will suffice to show that, for each projective variety,
\begin{equation}
\label{vanish}
\sigma(\cf_X) = 0
\end{equation}
in the homology group of $X$.
\par
If $X$ is a projective space, then
by hypothesis the top-dimensional component
of $\sigma(\cf_X)$ vanishes. From
Theorem~\ref{proj-re} we deduce equation~\ref{vanish}.
\par
In general we use induction on the dimension of $X$,
together with the Galois covering techniques of section~\ref{gct}.
Let $n$ be the dimension of $X$;
suppose that $\sigma(\cf_W) = 0$ for all varieties
of smaller dimension.
By Lemma~\ref{gct1},
we can construct a normal projective variety $Z$,
a finite morphism $\pi: Z \to X$,
and a finite morphism $\rho: X \to \bold{P}^n$,
such that the composition $\gamma: Z \to \bold{P}^n$ is Galois.
By Lemma~\ref{gct2}, $\bold{P}^n$
is the quotient of $Z$ by the Galois group $G$.
The characteristic function of $Z$ is fixed by the action
of the group; hence by naturality it is an element of
the invariant homology group.
By naturality and the inductive hypothesis
\begin{equation*}
\gamma_*\sigma(\cf_Z)
= \sigma\left(|G| \cdot \cf_{\bold{P}^n}
+ \text{ function supported on varieties of smaller dimension}
\right) = 0.
\end{equation*}
By Lemma~\ref{gct3},
$\gamma_*$ maps the invariant homology isomorphically
to the homology of $\bold{P}^n$. Hence
$\sigma(\cf_Z)$ must be zero.
Let $d$ be the degree of $\pi$.
Then
\begin{equation*}
0 = \pi_* \sigma(\cf_Z) =
\sigma\left(d \cdot \cf_X \right)
+ \sigma(\text{function supported on varieties of smaller dimension}).
\end{equation*}
By the inductive hypothesis the second term is zero.
Hence $\sigma(\cf_X) = 0$.
\end{pf*}
\par
Finally we wish to remark that the theorem and corollary are also
valid if we interpret all functors as emanating from the larger category
of compact complex algebraic varieties. Indeed, by
Chow's lemma \cite[p.\ 282]{shaf} an arbitrary compact variety is the image
of a projective variety via a birational morphism.
Thus we may prove an extended version of Theorem~\ref{main}
by an induction on dimension.
%
%


\begin{thebibliography}{9}

\bibitem{bs}
J. P. Brasselet and M. H. Schwartz,
{\em Sur les classes de Chern d'une ensemble analytique complexe},
Ast\'erisque 82--83 (1981), 93--148.

\bibitem{fulton}
W. Fulton,
{\em Rational equivalence on singular varieties},
Publ.\ Math.\ I. H. E. S. {\bf 45} (1975), 147--167.

\bibitem{int}
W. Fulton,
{\em Intersection Theory},
Ergebnisse der Math., 3 Ser.\ Vol.\ 2,
Springer, New York, 1984.

\bibitem{hardt}
R. M. Hardt,
{\em Triangulation of subanalytic sets and proper light
subanalytic maps}, Invent.\ Math.\ {\bf 38} (1977), 207--217.

\bibitem{harris}
J. Harris,
{\em Algebraic Geometry:  A First Course},
Springer, New York, 1992.

\bibitem{macp}
R. MacPherson,
{\em Chern classes for singular algebraic varieties},
Ann.\ Math.\ {\bf 4} (1974), 423--432.

\bibitem{mumford}
D. Mumford,
{\em The Red Book of Varieties and Schemes},
Lecture Notes in Math.\ 1358,
Springer, New York, 1988.

\bibitem{shaf}
I. Shafarevich,
{\em Basic Algebraic Geometry},
Springer, New York, 1977.

\bibitem{stewart}
I. Stewart,
{\em Galois Theory},
Chapman and Hall, New York, 1989.

\end{thebibliography}
\end{document}